\documentstyle[11pt,a4]{article}
\expandafter\chardef\csname pre amssym.def at\endcsname=\the\catcode`\@
\catcode`\@=11

\def\newsymbol#1#2#3#4#5{\let\next@\relax
 \ifnum#2=\@ne\let\next@\msafam@\else
 \ifnum#2=\tw@\let\next@\msbfam@\fi\fi
 \mathchardef#1="#3\next@#4#5}
 \def\hexnumber@#1{\ifcase#1 0\or 1\or 2\or 3\or 4\or 5\or 6\or 7\or 8\or
 9\or A\or B\or C\or D\or E\or F\fi}

\font\tenmsb=msbm10 scaled \magstep 1
\font\sevenmsb=msbm7 scaled \magstep 1
\font\fivemsb=msbm5 scaled \magstep 1
\newfam\msbfam
\textfont\msbfam=\tenmsb
\scriptfont\msbfam=\sevenmsb
\scriptscriptfont\msbfam=\fivemsb
\edef\msbfam@{\hexnumber@\msbfam}
\def\Bbb#1{\fam\msbfam\relax#1}

\def\widehat#1{\setbox\z@\hbox{$\m@th#1$}%
 \ifdim\wd\z@>\tw@ em\mathaccent"0\msbfam@5B{#1}%
\else\mathaccent"0362{#1}\fi} 

\font\teneufm=eufm10 scaled \magstep 1
\font\seveneufm=eufm7 scaled \magstep 1
\font\fiveeufm=eufm5 scaled \magstep 1
\newfam\eufmfam
\textfont\eufmfam=\teneufm
\scriptfont\eufmfam=\seveneufm
\scriptscriptfont\eufmfam=\fiveeufm
\def\frak#1{{\fam\eufmfam\relax#1}}

\newsymbol\subsetneqq 2324
\newsymbol\supsetneqq 2325

\catcode`\@=\csname pre amssym.def at\endcsname
 
\newtheorem{theorem}{Theorem}
\newtheorem{proposition}{Proposition}
 
\title{\bf HIDDEN SYMMETRY OF THE DIFFERENTIAL CALCULUS ON THE QUANTUM
MATRIX  SPACE}

\author{\sc S. Sinel'shchikov\thanks{Partially supported by the ISF grant
U2B200}\and \sc L. Vaksman\thanks{Partially supported by the AMS fSU grant
and by ISF grant U21200}}

\date{\it Mathematics Department, Institute for Low Temperature Physics and
Engineering,  47 Lenin Avenue, 310164 Kharkov, Ukraine}

\begin{document}
\begin{large}

\maketitle

 {\bf 1.}  This work solves a problem whose simple special case occurs in a
construction of a quantum unit ball of ${\Bbb C}^n$ (in the spirit of \cite{R}).
Within the framework of that theory, the automorphism group of the ball 
$SU(n,1)\subset SL(n+1)$ is essential. The problem is that the Wess-Zumino
differential calculus in quantum ${\Bbb C}^n$ \cite{WZ} at a first glance seems
to be only $U_q\frak{sl}_n $-invariant. In that particular case the lost 
$U_q\frak{sl}_{m+n}$-symmetry can be easily detected. The main result
of this work is in disclosing the hidden $U_q\frak{sl}_n $-symmetry for
bicovariant differential calculus in the quantum matrix space
$\hbox{Mat}(m,n)$. (Note that for $n=1$ we have the case of a ball).

 The aurhors are grateful to V. Akulov and G. Maltsiniotis for a helpful
discussion of the results.

\bigskip

 {\bf 2.} We start with recalling the definition of the Hopf algebra
$U_q\frak{sl}_N$, $N>1$, over the field ${\Bbb C}(q)$ of rational functions of an indeterminate $q$ \cite{D}, \cite{J}. (We follow the notations of \cite{CK}).

 For $i,j \in \{1,\ldots,N-1\}$ let

$$a_{ij}=\left\{\begin{array}{rl}{2,} & {i-j=0}\cr {-1,} &{\vert 
i-j\vert=1}\cr {0,} &{\vert i-j\vert >1.}\end{array}\right.$$

\smallskip

 The algebra $U_q\frak{sl}_N$ is defined by the generators
$\{E_i,\:F_i,\:K_i,\:K_i^{-1}\}$ and the relations

$$K_iK_j\,=\,K_jK_i,\quad K_iK_i^{-1}\,=\,K_i^{-1}K_i\,=\,1$$
$$K_iE_j\,=\,q^{a_{ij}}E_jK_i,\quad K_iF_j\,=\,q^{-a_{ij}}F_jK_i$$
$$E_iF_j\,-\,F_jE_i\:=\:\delta_{ij}(K_i-K_i^{-1})/(q-q^{-1})$$
$$E_i^2E_j\,-\,(q+q^{-1})E_iE_jE_i \,+\,E_jE_i^2\:=\:0, \quad \vert 
i-j\vert \,=\,1$$
$$F_i^2F_j\,-\,(q+q^{-1})F_iF_jF_i \,+\,F_jF_i^2\:=\:0, \quad \vert
i-j\vert \,=\,1$$
$$[E_i,E_j]\,=\,[F_i,F_j]\,=\,0, \quad \vert i-j\vert\,\ne\,1.$$

\smallskip

 A comultiplication $\Delta$, an antipode $S$ and a counit $\varepsilon$
are defined by
 
$$\Delta E_i\:=\:E_i\otimes 1\,+\,K_i\otimes E_i,\quad
\Delta F_i\:=\:F_i\otimes K_i^{-1}\,+\,1\otimes F_i,$$
$$\Delta K_i\,=\,K_i\otimes K_i,\quad S(E_i)\,=\,-K_i^{-1}E_i,$$
$$S(F_i)\,=\,-F_iK_i,\quad S(K_i)\,=\,K_i^{-1},$$
$$\varepsilon(E_i)=\varepsilon(F_i)=0,\quad\varepsilon(K_i)=1.$$

\bigskip

 {\bf 3.} Remind a description of a differential algebra
$\Omega^*(\hbox{Mat}(m,n))_q$ on a quantum matrix space \cite{CP} \cite{M1}.

 Let $i,j,i',j'\in\{1,2,\ldots,m+n\}$, and
 
$$\check{R}_{ij}^{i'j'}\;=\;\left\{\begin{array}{rl}{q^{-1},}&{i=j=i'=j'}\cr{1,}&
{i'=j\;\hbox{and}\;j'=i\;\hbox{and}\;i\ne j}\cr{q^{-1}-q,}&{i=i'\;\hbox{and}\;j=j'\;
\hbox{and}\;i<j}\cr{0,}&{\hbox{otherwise}}\end{array}\right.$$

 $\Omega^*(\hbox{Mat}(m,n))_q$ is given by the generators $\{t_a^\alpha\}$ and
relations

$$\sum_{\gamma,\delta}\check{R}_{\gamma\delta}^{\alpha\beta}t_a^\gamma
t_b^\delta\,=\,\sum_{c,d}\check{R}_{ab}^{cd}t_d^\beta t_c^\alpha$$
$$\sum_{a',b',\gamma',\delta'}\check{R}_{\gamma'\delta'}^{\alpha\beta}
\check{R}_{ab}^{a'b'}t_{a'}^{\gamma'}dt_{b'}^{\delta'}\;=\;dt_a^\alpha t_b^\beta$$
$$\sum_{a',b',\gamma',\delta'}\check{R}_{\gamma'\delta'}^{\alpha\beta}
\check{R}_{ab}^{a'b'}dt_{a'}^{\gamma'}dt_{b'}^{\delta'}\;=\;-dt_a^\alpha
dt_b^\beta$$

\noindent $(a,b,c,d,a',b'\,\in\{1,\ldots,n\};\quad
\alpha,\beta,\gamma,\delta,\gamma',\delta'\,\in\{1,\ldots,m\}).$

  Let us define a grading by
$\hbox{deg}(t_a^\alpha)\,=\,0,\;\hbox{deg}(dt_a^\alpha)\,=\,1$. With that,
${\Bbb C}[\hbox{Mat}(m,n)]_q\,=\,\Omega^0(\hbox{Mat}(m,n)))_q$ will stand 
for a subalgebra of zero degree elements .

\bigskip

  {\bf 4.} Let $A$ be a Hopf algebra and $F$ an algebra with unit and an
$A$-module the same time. $F$ is said to be a  $A$-module algebra \cite{A}
if the multiplication $m:\:F\otimes F\,\to\,F$ is a morphism of $A$-modules,
and $1\in F$ is an invariant (that is\\
$a(f_1f_2)\,=\,\sum\limits_j a_j^{'}f_1\otimes a_j^{''}f_2,\quad  a1\,=\,\varepsilon(a)1$ for all
$a\in A;\;f_1,f_2\in F$, with $\Delta(a)\,=\,\sum\limits_ja_j^{'}
\otimes a_j^{''}$).

  An important example of an $A$-module algebra appears if one supplies
$A^*$ with the structure of an $A$-module: $\langle af,b\rangle\,=\,\langle
f,ba\rangle,\;a,b\in A,\:f\in A^*$.

\bigskip

 {\bf 5.} Our immediate goal is to furnish ${\Bbb C}[\hbox{Mat}(m,n)]_q$ with
a structure of a $U_q\frak{sl}_{m+n}$-module algebra via an embedding ${\Bbb
C}[\hbox{Mat}(m,n)]_q\hookrightarrow (U_q\frak{sl}_{m+n})^*$.

 Let $\{e_{ij}\}$ be a standard basis in $\hbox{Mat}(m+n)$ and $\{f_{ij}\}$ 
the dual basis in $\hbox{Mat}(m+n)^*$. Consider a natural representation $\pi$ of
$U_q\frak{sl}_{m+n}$ :
$$\pi(E_i)\,=\,e_{i\:i+1},\quad \pi(F_i)\,=\,e_{i+1\:i},
\quad\pi(K_i)\:=\:qe_{ii}\,+\,q^{-1}e_{i+1\:i+1}\,+\,\sum_{j\ne i,i+1}e_{jj}.$$

 The matrix elements $u_{ij}\:=\:f_{ij}\pi\in(U_q\frak{sl}_{m+n})^*$ of
the natural representation may be treated as "coordinates" on the quantum group
$SL_{m+n}$ \cite{D}. To construct "coordinate" functions on a big cell of
the Grassmann manifold, we need the following elements of ${\Bbb
C}[\hbox{Mat}(m,n)]_q$
$$x(j_1,j_2,\ldots,j_m)\:=\:\sum\limits_{w\in
S_m}\,(-q)^{l(w)}u_{1j_{w(1)}}u_{2j_{w(2)}}\ldots u_{mj_{w(m)}},$$
with $1\le j_1<j_2<\ldots<j_m\le m+n$, and
$l(w)\;=\;\hbox{card}\{(a,b)\vert\:a<b\;\hbox{and}\;w(a)>w(b)\}$ being
the "length" of a permutation $w\in S_m$.

\medskip

  \begin{proposition}\hspace{-.5em}. $x(1,2,\ldots,m)$ is invertible in
$(U_q\frak{sl}_{m+n})^*$, and the map
$$t_a^\alpha\mapsto
x(1,2,\ldots,m)^{-1}x(1,\ldots,\widehat{m+1-\alpha},\ldots,m,m+a)$$
can be extended up to an embedding
$$i:\;{\Bbb C}[\hbox{Mat}(m,n)]_q\hookrightarrow(U_q\frak{sl}_{m+n})^*.$$
(The sign $\widehat{\ }$ here indicates the item in a list that should be
omitted).
\end{proposition}

\medskip

  Proposition 1 allows one to equip ${\Bbb C}[\hbox{Mat}(m,n)]_q$ with
the structure of a $U_q\frak{sl}_{m+n}$-module algebra :\\
$$i\xi t_a^\alpha\,=\,\xi it_a^\alpha, \qquad \xi\in U_q\frak{sl}_{m+n},\;
a\in\{1,\ldots,n\},\;\alpha\in\{1,\ldots,m\}.$$

\bigskip

  {\bf 6.} The main result of our work is the following
 \begin{theorem}\hspace{-.5em}. $\Omega^*(\hbox{Mat}(m,n))_q$ admits a
unique structure of a $U_q\frak{sl}_{m+n}$-module algebra such that
the embedding
$$i:\:{\Bbb C}[\hbox{Mat}(m,n)]_q\hookrightarrow
\Omega^*(\hbox{Mat}(m,n))_q$$
and the differential
$$d:\;\Omega^*(\hbox{Mat}(m,n))_q\,\to\,\Omega^*(\hbox{Mat}(m,n))_q$$
are the morphisms of $U_q\frak{sl}_{m+n}$-modules.
\end{theorem}

\medskip

  {\sc Remark 1.} The bicovariance of the differential calculus on the
quantum matrix space allows one to equip the algebra
$\Omega^*(\hbox{Mat}(m,n))_q$ with a structure of
$U_q\frak{s}(\frak{gl}_m\times\frak{gl}_n)$-module, which is compatible
with multiplication in $\Omega^*(\hbox{Mat}(m,n))_q$ and differential $d$.
Theorem 1 implies that $\Omega^*(\hbox{Mat}(m,n))_q$ possess an additional
hidden symmetry since
$U_q\frak{sl}_{m+n}\,\supsetneqq\,U_q\frak{s}(\frak{gl}_m
\times\frak{gl}_n)$. 

\medskip

 {\sc Remark 2.} Let $q_0\in{\Bbb C}$ and $q_0$ is not a root of unity.
It follows from the explicit formulae for $E_mt_a^\alpha,\;F_mt_a^\alpha,
\;K_m^{\pm 1}t_a^\alpha,\quad a\in\{1,\ldots,n\},\quad\alpha\in\{1,
\ldots,m\}$, that the "specialization" $\Omega^*(\hbox{Mat}(m,n))_{q_0}$ is a
$U_{q_0}\frak{sl}_{m+n}$-module algebra.

\bigskip

 {\bf 7.} Supply the algebra $U_q\frak{sl}_{m+n}$ with a grading as follows:
$$\hbox{deg}(K_i)\,=\,\hbox{deg}(E_i)\,=\,\hbox{deg}(F_i)\,=\,0,\quad
\hbox{for}\quad i\ne m,$$
$$\hbox{deg}(K_m)\,=\,0,\quad \hbox{deg}(E_m)\,=\,1,\quad
\hbox{deg}(F_m)\,=\,0.$$

  The proofs of Proposition 1 and Theorem 1 reduce to the construction of
graded $U_q\frak{sl}_{m+n}$-modules which are dual respectively to the
modules of functions $\Omega^0(\hbox{Mat}(m,n))_q$ and that of 1-forms
$\Omega^1(\hbox{Mat}(m,n))_q$. The dual modules are defined by their
generators and correlations. While proving the completeness of the
correlation list we implement the "limit specialization" $q_0=1$
(see \cite{CK}, p. 476).

 The passage from the order one differential calculus
$\Omega^0(\hbox{Mat}(m,n))_q\,\stackrel{d}{\to}\,
\Omega^1(\hbox{Mat}(m,n))_q$ to $\Omega^*(\hbox{Mat}(m,n))_q$ is done
via a universal argument described in a paper by G. Maltsiniotis \cite{M2}.
This argument doesn't break $U_q\frak{sl}_{m+n}$-symmetry.

\bigskip

 {\bf 8.} Our approach to the construction of order one differential
calculus is completely analogous to that of V. Drinfeld \cite{D} used
initially to produce the algebra of functions on a quantum group by means
of a universal enveloping algebra.

\bigskip

  {\bf 9.} The space of matrices is the simplest example of an irreducible
prehomogeneous vector space of parabolic type \cite{Ki}. Such space can
be also associated to a pair constituted by a Dynkin diagram of a simple
Lie algebra ${\cal G}$ and a distinguished vertex of this diagram. Our
method can work as an efficient tool for producing $U_q{\cal G}$-invariant
differential calculi on the above prehomogeneous vector spaces.

 Note that $U_q{\cal G}$-module algebras of polynomials on quantum
prehomogeneous spaces of parabolic type were considered in a recent
work of M. S. Kebe \cite{Ke}.

\bigskip

{\bf Acknowledgement.}  The  authors would like to express their gratitude
to Prof. A. Boutet de Monvel at  University Paris VII  for the warm hospitality
during the work on this paper.

\begin{center}

\end{center}

\end{large}

\begin{thebibliography}{99}

\bibitem{A} E. Abe. Hopf Algebras, Cambridge Tracts in Mathematics 74,
Cambridge University Press, Cambridge, 1980. 

\bibitem{CP} V. Chari and A. Pressley. A guide to quantum groups. Cambridge
Univ. Press, 1994.

\bibitem{CK} C. De Concini and V. G. Kac, {\it Representations of quantum groups
at roots of $1$}, in Operator Algebras, Unitary representations, Enveloping
Algebras and Invariant Theory, A. Connes, M. Duf\/lo, A. Joseph, R. Rentschler
(eds.), pp. 471- 506, Birkh\"{a}user, Boston, 1990.

\bibitem{D} V. G. Drinfel'd , {\it Quantum groups}, in Proceedings of
International Congress of Mathematicians, Berkeley, 1986, A. M. Gleason
(ed.), pp. 798 - 820, American Mathematical Society, Providence, RI, 1987.

\bibitem{J} M. Jimbo, {\it Quantum $R$-matrix related to the generalized
Toda system: an algebraic approach}, in Field Theory, Quantum Gravity and
Strings, H. ~J. de Vega \& N. Sanchez (eds.), Lecture Notes in Phys., {\bf 246},
335 - 361, Springer, Berlin,1986.

\bibitem{Ke} M\^alek Stefan Kebe, $\cal{O}$-{\it alg\`ebres quantiques},
C. R. Acad. Sci. Paris, {\bf 322} (1996), S\'erie 1, 1 - 4.

\bibitem{Ki} T. Kimura, {\it A classification theory of Prehomogeneous
vector spaces}, Advanced Studies in Pure Mathematics {\bf 14}, 1988,
Representations of Lie groups, Kyoto, Hiroshima, 1986, p. 223 - 256.

\bibitem{M1} G. Maltsiniotis, {\it Groupes quantiques et structures
diff\'erentielles}, C. R. Acad. Sci., Paris, S\'erie I, {\bf 311} (1990), 831 - 834.

\bibitem{M2} G. Maltsiniotis, {\it Le langage des espaces et des groupes
quantiques}, Comm. Math. Phys., {\bf 151} (1993), 275 - 302.

\bibitem{R} W. Rudin. Function Theory in the Unit Ball of ${\Bbb C}^n$.
Springer -- Verlag, N.-Y., 1980.

\bibitem{WZ} J. Wess and B. Zumino, {\it Covariant differential calculus on the quantum hyperplane}, Nucl. Phys. B, Proc. Suppl. {\bf 18 B}, 1991, 302 - 312.

\end{thebibliography}
\end{document}